\input amstex
\input Amstex-document.sty

\pageno 711

\def\shortauthor{N. Lerner}
\def\shorttitle{Solvability}
\def\re{\operatorname{Re}}
\def\im{\operatorname{Im}}
\def\supp{\operatorname{supp}}
\font\boldgreeknormal = cmmib10 scaled 1000
\font\boldgreek = cmmib10 scaled 1440

\topmatter

\title\nofrills{\boldHuge Solving Pseudo-Differential Equations}
\endtitle

\author \Large Nicolas Lerner* \endauthor

\thanks *University of Rennes, Universit\'e de Rennes 1, Irmar, Campus de Beaulieu, 35042 Rennes
cedex, France. E-mail: lerner\@univ-rennes1.fr \endthanks

\abstract\nofrills \centerline{\boldnormal Abstract}

\vskip 4.5mm

{\ninepoint In 1957, Hans Lewy constructed a counterexample
showing that very simple and natural differential equations can
fail to have local solutions. A geometric interpretation and a
generalization of this counterexample were given in 1960 by
L.H\"ormander. In the early seventies, L.Nirenberg and F.Treves
proposed a geometric condition on the principal symbol, the
so-called condition $(\psi)$, and  provided strong arguments
suggesting that it should be equivalent to local solvability. The
necessity of condition $(\psi)$ for solvability of
pseudo-differential equations was proved  by L.H\"ormander in
1981. The sufficiency of condition $(\psi)$ for solvability of
differential equations was proved by R.Beals and C.Fefferman in
1973. For differential equations in any dimension and for
pseudo-differential equations in two dimensions, it was shown more
precisely that $(\psi)$ implies solvability with a loss of one
derivative with respect to the elliptic case: for instance, for a
complex vector field $X$ satisfying $(\psi)$, $f\in
L^2_{\text{loc}}$, the equation $Xu=f$ has a solution $u\in
L^2_{\text{loc}}$.

In 1994, it was  proved by N.L. that condition $(\psi)$ does not
imply solvability with loss of one derivative for
pseudo-differential equations, contradicting repeated claims by
several authors. However in 1996, N.Dencker proved that  these
counterexamples were indeed solvable, but with a loss of two
derivatives. We shall explore the structure of this phenomenon
from both sides: on the one hand, there are first-order
pseudo-differential equations satisfying condition $(\psi)$ such
that no $L^2_{\text{loc}}$ solution can be found with some source
in $L^2_{\text{loc}}$. On the other hand, we shall see that, for
these examples, there exists a solution in the Sobolev space
$H^{-1}_{\text{loc}}$.

The sufficiency of condition $(\psi)$ for solvability of
pseudo-differential equations in three or more dimensions is still
an open problem. In 2001, N.Dencker announced that he has proved
that condition $(\psi)$ implies solvability (with a loss of two
derivatives), settling the Nirenberg-Treves conjecture. Although
his paper contains several bright and new ideas, it is the opinion
of the author of these lines that a number of points in his
article need clarification.

\vskip 4.5mm


}
\endabstract
\endtopmatter

\document

\specialhead \noindent \boldLARGE 1. From  Hans Lewy to Nirenberg-Treves' \\ \text{ } \hskip 4mm condition
({\boldgreek \char32}) \endspecialhead

\subhead Year 1957 \endsubhead

The Hans Lewy operator $L_0$, introduced in \cite{20}, is the
following  complex vector field in ${\Bbb R}^3$
$$
L_0=\frac{\partial} {{\partial} x_1}+i\frac{{\partial}}{{\partial}
x_2}+i(x_1+ix_2)\frac{{\partial}}{{\partial} x_3}. \tag 1.1$$
There exists $f\in C^\infty$ such that the equation
$$
L_0u=f \tag 1.2$$ has no distribution solution, even locally. This
discovery came as a great shock for several reasons. First of all,
$L_0$ has a very simple expression
 and  is natural as
the Cauchy-Riemann operator on the boundary of the pseudo-convex
domain
$$
\{(z_1,z_2)\in{\Bbb C}^2,\vert{z_1}\vert^2+2\im{z_2}<0\}.
$$
Moreover $L_0$ is a non-vanishing vector field so that no
pathological behaviour related to multiple characteristics is to
be expected. In the fifties, it was certainly the conventional
wisdom that any ``reasonable" operator should be locally solvable,
and obviously (1.1) was indeed very reasonable so the conclusion
was that, once more, the {CW} should be revisited. One of the
questions posed by  such a counterexample was to find some
geometric explanation for this phenomenon.

\subhead 1960 \endsubhead

This was done in 1960 by L.H\"ormander in \cite{7} who proved that
if $p$ is the symbol of a differential operator such that, at some
point $(x,\xi)$ in the  cotangent bundle,
$$
p(x,\xi)=0\quad \text{and}\quad \{\re p,\im p\}(x,\xi) >0, \tag
1.3$$ then the operator $P$ with principal symbol $p$ is not
locally solvable at  $x$; in fact, there exists $f\in C^\infty$
such that, for any neighborhood $V$ of $x$ the equation $Pu=f$ has
no  solution $u\in\Cal D' (V)$. Of course, in the case of
differential operators, the sign $>0$ in (1.3) can be replaced by
$\not= 0$ since the Poisson bracket $\{\re p,\im p\}$ is then an
homogeneous polynomial with odd degree in the variable $\xi$.
Nevertheless, it appeared later (in \cite{8}) that the same
statement is true for pseudo-differential operators, so we keep it
that way. Since the symbol of $-iL_0$ is
$\xi_1-x_2\xi_3+i(\xi_2+x_1\xi_3)$, and the Poisson bracket $
\{\xi_1-x_2\xi_3,\xi_2+x_1\xi_3\}=2\xi_3, $ the assumption (1.3)
is fulfilled for $L_0$ at any point $x$ in the base and the
nonsolvability property follows. This gives a necessary condition
for local solvability of pseudo-differential equations: a locally
solvable operator $P$ with principal symbol $p$ should satisfy
$$
 \{\re p,\im p\}(x,\xi) \le 0\quad
\text{at} \quad p(x,\xi)=0. \tag 1.4$$ Naturally, condition (1.4)
is far from being sufficient for solvability (see e.g. the
nonsolvable $M_3$ below in (1.5)). After the papers \cite{20},
\cite{7}, the curiosity of the mathematical community was aroused
in search of a geometric condition on the principal symbol,
characterizing local solvability of principal type operators. It
is important to note that for principal type operators with a real
principal symbol, such as a non-vanishing real vector field, or
the wave equation, local solvability was known after the 1955
paper of L.H\"ormander in \cite{6}. In fact these results extend
quite easily to the pseudo-differential real principal type case.
As shown by the Hans Lewy counterexample and the necessary
condition (1.4), the matters are quite different for
complex-valued symbols.

\subhead 1963 \endsubhead

It is certainly helpful to look now  at some simple models. For
$t,x\in{\Bbb R}$, with the usual notations
$$
D_t=-i{\partial}_t, \quad\widehat{(\vert D_x\vert
u)}(\xi)=\vert{\xi}\vert\hat{u}(\xi),
$$
where $\hat u$ is the $x$-Fourier transform of $u$, $l\in {\Bbb
N}$, let us consider the operators defined by
$$
M_l=D_t+it^{l}D_x,\quad N_l=D_t+it^{l}\vert{D_x}\vert. \tag 1.5$$
It is indeed rather easy to prove that, for $k\in{\Bbb N}$,
$M_{2k},N_{2k}, N_{2k+1}^*$ are solvable whereas $M_{2k+1},
N_{2k+1}$ are nonsolvable. In particular, the operators $M_1,N_1$
satisfy $(1.3)$. On the other hand, the operator
$N_1^*=D_t-it\vert{D_x}\vert$ is indeed solvable since its adjoint
operator $N_1$ verifies the a priori estimate
$$
T\Vert N_1u\Vert_{L^2({\Bbb R}^2)}\ge \Vert{u}\Vert_{L^2({\Bbb
R}^2)},
$$
for a smooth compactly supported $u$ vanishing for $\vert{t}\vert
\ge T/2$. No such estimate is satisfied by $N_1^*u$ since its
$x$-Fourier transform is
$$
-i{\partial}_t v-it\vert{\xi}\vert v=(-i) ({\partial}_t
v+t\vert\xi\vert v),
$$
where $v$ is the $x$-Fourier transform of $u$. A solution of
$N_1^*u=0$ is thus given by the inverse Fourier transform of $
e^{-t^2\vert{\xi}\vert/2} $, ruining solvability for the operator
$N_1$. A complete study of solvability  properties of the models
$M_l$ was done in \cite{23} by  L.Nirenberg and F.Treves, who also
provided a sufficient condition of solvability for vector fields;
the analytic-hypoellipticity properties of these operators were
also studied in a paper by S.Mizohata \cite{21}.

\subhead 1971 \endsubhead

The ODE-like examples (1.5) led L.Nirenberg and F.Treves in
\cite{24--25--26} to formulate a conjecture and to prove it in a
number of cases, providing strong grounds in its favour. To
explain this, let us look simply at the operator
$$
L=D_t+i q(t,x,D_x), \tag 1.6$$ where $q$ is a real-valued
first-order symbol. The symbol of $L$ is thus $\tau+iq(t,x,\xi)$.
The bicharacteristic curves of the real part are oriented straight
lines with direction ${\partial}/{\partial} t$; now we examine the
variations of the imaginary part $q(t,x,\xi)$ along these lines.
It amounts only to check the functions $t\mapsto q(t,x,\xi)$ for
fixed $(x,\xi)$. The good cases in (1.5) (when solvability holds)
are $t^{2k}\xi,-t^{2k+1}\vert{\xi}\vert$: when $t$ increases these
functions do not change sign from $-$ to $+$. The bad cases are
$t^{2k+1}\vert{\xi}\vert$: when $t$ increases these functions do
change sign from $-$ to $+$; in particular, the nonsolvable case
(1.3), tackled in \cite{8}, corresponds to a change of sign of
$\im{p}$ from $-$ to $+$ at a simple zero. The general formulation
of condition $(\psi)$ for a principal type operator with principal
symbol $p$ is as follows: for all $z\in{\Bbb C}$, $\im (zp)$ does
not change sign from $-$ to $+$ along the oriented
bicharacteristic curves of $\re (zp)$. It is a remarkable and
non-trivial fact that this condition is invariant by
multiplication by an elliptic factor as well as by composition
with an homogeneous canonical transformation. The {\it
Nirenberg-Treves conjecture}, proved in several cases in
\cite{24--25--26}, such as for differential operators with
analytic coefficients, states that, {\it for a principal type
pseudo-differential equation, condition $(\psi)$ is equivalent to
local solvability.}
\par
The paper \cite{25} introduced a radically new method of proof of
energy estimates for the adjoint operator $L^*$ based on a
factorization of $q$ in (1.6): whenever
$$
q(t,x,\xi)=a(t,x,\xi)b(x,\xi) \tag 1.7$$ with $a\le 0$ of order 0
and $b$ of order 1, then the operator $L$ in (1.6) is locally
solvable. Looking simply at the ODE
$$
D_t+i
a(t,x,\xi)b(x,\xi)=(-i)\bigl({\partial}_t-a(t,x,\xi)b(x,\xi)\bigr),
\tag 1.8$$ it is clear that in the region $\{b(x,\xi)\ge 0\}$, the
forward Cauchy problem for (1.8) is well posed, whereas in
$\{b(x,\xi)\le 0\}$, well-posedness holds for the backward Cauchy
problem. This remark led L.Nirenberg and F.Treves to use  as a
multiplier in the energy method the sign of the operator with
symbol $b$. They were also able to provide the proper commutator
estimates to handle the remainder terms generated by this
operator-theoretic method. Although a factorization (1.7) can be
obtained for differential operators with analytic regularity
satisfying condition $(\psi)$, such a factorization is not true in
the $C^\infty$ case. Incidentally, one should note that for
differential operators, condition $(\psi)$ is equivalent to ruling
out any change of sign of $\im p$ along the bicharacteristics of
$\re p$ (the latter condition is called condition $(P)$); this
fact is due to the identity $ p(x,-\xi)=(-1)^mp(x,\xi), $ valid
for an homogeneous polynomial of degree $m$ in the variable $\xi$.
\par
Using the Malgrange-Weierstrass theorem on normal forms of
complex-valued non-degenerate $C^\infty$ functions and the Egorov
theorem on quantization of homogeneous canonical transformations,
there is no loss of generality considering only first order
operators of type $(1.6)$. The expression of condition $(\psi)$
for $L$ is then very simple since it reads
$$
q(t,x,\xi)<0\quad \text{and}\quad s>t \Longrightarrow
q(s,x,\xi)\le 0. \tag 1.9$$ Note that the expression of condition
$(P)$ for $L$ is simply $ q(t,x,\xi)q(s,x,\xi)\ge 0. $ Much later
in 1988, N.Lerner  \cite{14} proved the sufficiency of condition
$(\psi)$ for local solvability of pseudo-differential equations in
two dimensions and as well for the classical oblique-derivative
problem \cite{15}. The method of proof of these results is based
upon a factorization analogous to (1.7) but where $b(x,\xi)$ is
replaced by $\beta(t,x)\vert{\xi}\vert$ and $\beta$ is a smooth
function such that $t\mapsto \beta(t,x)$ does not change sign from
$+$ to $-$ when $t$ increases. Then a properly defined sign of
$\beta(t,x)$ appears as a non-decreasing operator and the
Nirenberg-Treves energy method can be adapted to this situation.

\subhead 1973 \endsubhead

At this date, R.Beals and C.Fefferman \cite{1} took as a starting
point the previous results of L.Nirenberg and F.Treves and,
removing the analyticity assumption, they were able to prove the
sufficiency of condition $(P)$ for local solvability, obtaining
thus the sufficiency of condition $(\psi)$ for local solvability
of differential equations. The  key ingredient was a drastically
new vision of the pseudo-differential calculus, defined to obtain
the factorization (1.7) in regions of the phase space much smaller
than cones or semi-classical ``boxes" $\{(x,\xi), \vert{x}\vert
\le 1,\vert{\xi}\vert\le h^{-1}\}$. Considering the family
$\bigl\{q(t,x,\xi)\bigr\}_{t\in[-1,1]}$ of classical homogeneous
symbols of order 1, they define, via a Calder\'on-Zygmund
decomposition, a pseudo-differential calculus depending on the
family $\{q(t,\cdot)\}$, in which all these symbols are first
order but also such that, at some level $t_0$, some ellipticity
property of $q(t_0,\cdot)$ or $\nabla_{x,\xi}q(t_0,\cdot)$ is
satisfied. Condition $(P)$ then implies easily a factorization of
type (1.7) and the Nirenberg-Treves energy method can be used. It
is interesting to notice that some versions of these new
pseudo-differential calculi were used later on for the proof of
the Fefferman-Phong inequality \cite{5}. In fact, the proof of
R.Beals and C.Fefferman marked the day when microlocal analysis
stopped being only homogeneous or semi-classical, thanks to
methods of harmonic analysis such as Calder\'on-Zygmund
decomposition made compatible with the Heisenberg uncertainty
principle.

\subhead 1978 \endsubhead

Going back to solvability problems, the existence of $C^\infty$
solutions for $C^\infty$ sources was proved by L.H\"ormander in
\cite{9} for pseudo-differential equations  satisfying condition
$(P)$. For such an operator $P$ of order $m$, satisfying also a
non-trapping condition, a semi-global existence theorem was
proved, with a loss of $1+\epsilon$ derivatives, with
$\epsilon>0$. Following an idea given by R.D.Moyer \cite{22} for a
result in two dimensions, L.H\"ormander proved in  \cite{10} that
condition $(\psi)$ is necessary for local solvability: assuming
that condition $(\psi)$ is not satisfied for a principal type
operator $P$, he was able to construct
 approximate non-trivial solutions $u$ for the adjoint equation
$P^*u=0$, which implies that $P$ is not solvable. Although the
construction is elementary for the model  operators $N_{2k+1}$ in
(1.5) (as sketched above for $N_1$ in our 1963 section), the
multidimensional proof is rather involved and based upon a
geometrical optics method adapted to the complex case. The details
can be found in the proof of theorem $26.4.7'$ of \cite{11}.
\par
We refer the reader to the paper \cite{13} for a more detailed
historical overview of this problem. On the other hand, it is
clear that our interest is focused on solvability in the
$C^\infty$ category. Let us nevertheless recall that the
sufficiency of condition $(\psi)$ in the analytic category (for
microdifferential operators acting on microfunctions) was proved
by J.-M.Tr\'epreau \cite{27} (see also \cite{12}, chapter
{\sevenrm VII}).

\specialhead \noindent \boldLARGE 2. Counting the loss of
derivatives \endspecialhead

\subhead {Condition ({\boldgreeknormal \char32}) does not imply solvability with loss of one derivative}
\endsubhead

Let us consider a principal-type pseudo-differential operator $L$
of order $m$. We shall say that $L$ is locally solvable with a
loss of $\mu$ derivatives whenever the equation $Lu=f$ has a local
solution $u$ in the Sobolev space $H^{s+m-\mu}$ for a source $f$
in $H^s$. Note that the loss is zero if and only if $L$ is
elliptic. Since for the simplest principal type equation
${\partial}/{\partial}x_1$, the loss of derivatives is 1, we shall
consider that 1 is the ``ordinary" loss of derivatives. When $L$
satisfies condition $(P)$ (e.g. if $L$ is a differential operator
satisfying condition $(\psi)$), or when $L$ satisfies condition
$(\psi)$ in two dimensions, the estimates
$$
C\Vert{L^* u}\Vert_{H^s}\ge \Vert{u}\Vert_{H^{s+m-1}}, \tag 2.1$$
valid for smooth compactly supported $u$ with a small enough
support, imply local solvability with loss of 1 derivative, the
ordinary loss referred to above. For many years, repeated claims
were made that condition $(\psi)$ for $L$ implies (2.1), that is
solvability with loss of 1 derivative. It turned out that these
claims were wrong, as shown in \cite{16} by the following result
(see also section 6 in the survey \cite{13}). \proclaim{Theorem
2.1} There exists a principal type first-order pseudo-differential
operator $L$ in three dimensions, satisfying condition $(\psi)$, a
sequence $u_k$ of $C_{c}^\infty$ functions with $
\supp{u_k}\subset\{x\in{\Bbb R}^3,\vert{x}\vert\le 1/k\} $ such
that
$$
\Vert{u_k}\Vert_{L^2({\Bbb R}^3)} =1,\quad \lim_{k\rightarrow
+\infty}\Vert{L^* u_k}\Vert_{L^2({\Bbb R}^3)}=0. \tag 2.2$$
\endproclaim
As a consequence, for this $L$, there exists $f\in L^2$ such that
the equation $Lu=f$ has no local solution $u$ in $L^2$. We shall
now briefly examine some of the main features of this
counterexample, leaving aside the technicalities which can be
found in the papers quoted above. Let us try, with
$(t,x,y)\in{\Bbb R}^3$,
$$
L=D_t-ia(t)\bigl(D_x+ H(t)V(x)\vert{D_y}\vert\bigr), \tag 2.3$$
with $H=\bold{1}_{{\Bbb R}_+}$, $C^\infty({\Bbb R})\ni V\ge 0$,
$C^\infty({\Bbb R})\ni a\ge 0$ flat at $0$. Since the function
$q(t,x,y,\xi,\eta)=-a(t)\bigl(\xi+ H(t) V(x)\vert\eta\vert\bigr)$
satisfies (1.9) as the product of the non-positive function
$-a(t)$ by the non-decreasing function $t\mapsto \xi+
H(t)V(x)\vert\eta\vert$, the operator $L$ satisfies condition
$(\psi)$. To simplify the exposition, let us assume that $a\equiv
1$, which introduces a rather unimportant singularity in the
$t$-variable, let us replace $\vert{D_y}\vert$ by a positive
(large) parameter $\Lambda$, which allows us to work now only with
the two real variables $t,x$ and let us set $W=\Lambda V$. We are
looking for a non-trivial solution $u(t,x)$ of $L^*u=0$, which
means then
$$
{\partial}_t u=\cases D_x u,&\text{for $t<0,$}
\\
\bigl(D_x+W(x)\bigr)u,&\text{for $t>0.$}
\endcases
$$
The operator $D_x+W$ is unitarily equivalent to $D_x$: with
$A'(x)=W(x)$, we have $ D_x+W(x)=e^{-iA(x)}D_xe^{iA(x)}, $ so that
the negative eigenspace of the operator $D_x+W(x)$ is $ \{v\in
L^2({\Bbb R}), \supp\widehat{\hskip1pt
e^{iA}v\hskip1pt}\subset{\Bbb R}_-\}. $ Since we want $u$ to decay
when $t\rightarrow\pm\infty$, we need to choose $v_1,v_2\in
L^2({\Bbb R})$, such that
$$
u(t,x)=\cases e^{tD_x} v_1,\quad\supp{\widehat{\hskip1pt
v_1}}\subset{\Bbb R}_+&\text{for $t<0,$}
\\
e^{t(D_x+W)}v_2, \quad\supp{\widehat{e^{iA}v_2}}\subset{\Bbb R}_-
&\text{for $t>0.$}
\endcases
\tag 2.4$$ We shall not be able to choose $v_1=v_2$ in (2.4), so
we could only hope for $L^* u$ to be small if
$\Vert{v_2-v_1}\Vert_{L^2({\Bbb R})}$ is small. Thus this
counterexample is likely to  work if the unit spheres of the
vector spaces
$$
E_1^+=\{v\in L^2({\Bbb R}),\supp\widehat v\subset{\Bbb R}_+\}
\quad\text{and}\quad E_2^-=\{v\in L^2({\Bbb R}),\supp\widehat
{e^{iA}v\hskip1pt} \subset{\Bbb R}_-\}
$$
are close. Note that since $W\ge 0$, we get  $E_1^+\cap
E_2^-=\{0\}$: in fact, with $L^2({\Bbb R})$ scalar products, we
have
$$
v\in E_1^+\cap E_2^-\Longrightarrow 0\overset{v\in E_1^+}\to\le
\langle Dv,v \rangle \overset{0\le W}\to\le \langle (D+W)v,v
\rangle \overset{v\in E_2^-}\to\le 0 \Longrightarrow \langle Dv,v
\rangle=0
$$
which gives $v=0$ since $v\in E_1^+$. Nevertheless, the ``angle"
between $E_1^+$ and $E_2^-$ could be small for a careful choice of
a positive $W$. It turns out that $W_0(x)=\pi\delta_0(x)$ is such
a choice. Of course, several problems remain such as regularize
$W_0$ in such a way that it becomes a first-order semi-classical
symbol, redo the same construction with a smooth function $a$ flat
at 0 and various other things.
\par
Anyhow, these difficulties eventually turn out to be only
technical, and {\it in fine}, the actual reason for which theorem
2.1 is true is simply that the positive eigenspace of $D_x$ (i.e.
$L^2({\Bbb R})$ functions whose Fourier transform is supported in
${\Bbb R}_+$) could be arbitrarily close to the negative
eigenspace of $D_x+W(x)$ for some non-negative $W$, triggering
nonsolvability in $L^2$ for the three-dimensional  model operator
$$
D_t-ia(t)\bigl(D_x+\bold{1}_{{\Bbb R}_+}(t)
W(x)\vert{D_y}\vert\bigr), \tag 2.5$$ where $a$ is some
non-negative function, flat at $0$. This phenomenon is called the
``drift" in \cite{16} and could not occur for differential
operators or for pseudo-differential operators in two dimensions.
A more geometric point of view is that for a principal type symbol
$p$, satisfying condition $(\psi)$, one may  have
bicharacteristics of $\re p$ which stay in the set $\{\im p=0\}$.
This can even occur for operators satisfying condition $(P)$.
However condition $(P)$ ensures that the nearby bicharacteristics
of $\re p$ stay either in $\{\im p\ge 0\}$ or in $\{\im p\le 0\}$.
This is no longer the case when condition $(\psi)$ holds, although
the bicharacteristics  are not allowed to pass from $\{\im p< 0\}$
to $\{\im p> 0\}$. The situation of having a bicharacteristic of
$\re p$ staying in $\{\im p=0\}$ will generically trigger the
drift phenomenon mentioned above when condition $(P)$ does not
hold. So the counterexamples to solvability with loss of one
derivative are in fact very close to operators satisfying
condition $(P)$.
\par
A related remark is that the ODE-like solvable models in (1.5) do
not catch the generality allowed by condition $(\psi)$. Even for
subelliptic operators, whose tranposed are of course locally
solvable, it is known that other model operators than
$M_{2k},N_{l}$ can occur. In particular the three-dimensional
models $ D_t+it^{2k}(D_x+t^{2l+1}x^{2m}\vert{D_y}\vert), $ where
$k,l,m$ are non-negative integers are indeed subelliptic and are
not reducible to (1.5) (see chapter 27 in \cite{11} and the remark
before corollary 27.2.4 there).

\subhead {Solvability with loss of two derivatives} \endsubhead

Although theorem 2.1 demonstrates that{~}condition $(\psi)$ does
not imply solvability with loss of one derivative, the
counterexamples constructed in this theorem are indeed solvable,
but with a loss of two derivatives, as proven by N.Dencker in 1996
\cite{2}. The same author gave a generalization of his results in
\cite{3} and later on, analogous results were given in \cite{17}.
\par
A measurable function $p(t,x,\xi)$ defined on ${\Bbb R}\times{\Bbb
R}^n\times {\Bbb R}^n$ will be called in the next theorem  a
symbol of order $m$ whenever, for all $(\alpha, \beta)\in{\Bbb
N}^n\times {\Bbb N}^n$
$$
\sup_{(t,x,\xi)\in {\Bbb R}\times{\Bbb R}^n\times {\Bbb R}^n }
\vert{({\partial}_x^\alpha{\partial}_\xi^\beta p)(t,x,\xi)}\vert
(1+\vert\xi\vert)^{-m+\vert\beta\vert}<+\infty. \tag 2.6$$
\proclaim{Theorem 2.2} Let $a(t,x,\xi)$ be a non-positive symbol
of order $0$, $b(t,x,\xi)$ be a real-valued symbol of order $1$
such that ${\partial}_t b\ge 0$, and $r(t,x,\xi)$ be a
(complex-valued) symbol of order $0$. Then the operator
$$
L=D_t+ia(t,x,D_x) b(t,x,D_x)+r(t,x,D_x) \tag 2.7$$ is locally
solvable with a loss of two derivatives. Since the counterexamples
constructed in theorem $2.1$ are in fact of type $(2.7)$, they are
locally  solvable with a loss of two derivatives.
\endproclaim
In fact, for all points in ${\Bbb R}^{n+1}$, there exists a
neighborhood $V$, a positive constant $C$ such that, for all $u\in
C_c^\infty(V)$
$$
C\Vert{L^* u}\Vert_{H^0} \ge \Vert{u}\Vert_{H^{-1}}. \tag 2.8$$
This estimate actually represents a loss of two derivatives for
the first-order $L$; the  estimate with loss of 0 derivative would
be $ \Vert{L^* u}\Vert_{H^0} \gtrsim \Vert{u}\Vert_{H^{1}}, $ the
estimate with loss of one derivative would be $ \Vert{L^*
u}\Vert_{H^0} \gtrsim \Vert{u}\Vert_{H^0}, $ and both are false,
the first because $L^*$ is not elliptic, the second from theorem
2.1. The proof of theorem 2.2 is essentially based upon the energy
method which boils down to compute for all $T\in{\Bbb R}$
$$
\re\langle L^* u, i B u+i H(t-T)u\rangle_{L^2({\Bbb R}^{n+1})}
$$
where $B=b(t,x,D_x)$. Some complications occur in the proof from
the rather weak assumption ${\partial}_tb\ge 0$ and also from the
lower order terms. Anyhow, the correct multiplier is essentially
given by $b(t,x,D_x)$. Theorem 2.2 can be proved for much more
general classes
 of pseudo-differential
operators than those given by (2.6). As a consequence, it can be
extended naturally to contain the solvability result under
condition $(P)$ (but with a loss of two derivatives, see e.g.
theorem 3.4 in \cite{17}).

\subhead {Miscellaneous results} \endsubhead

Let us  mention that the operator (1.6) is solvable with a loss of
one derivative (the ordinary loss) if condition $(\psi)$ is
satisfied (i.e. (1.9)) as well as the  extra condition
$$
\vert{{\partial}_x q(t,x,\xi)}\vert^2\vert{\xi}\vert^{-1} +
\vert{{\partial}_\xi q(t,x,\xi)}\vert^2\vert{\xi}\vert\le C
\vert{{\partial}_t q(t,x,\xi)}\vert\quad \text{when}\quad
q(t,x,\xi)=0.
$$
This result is proved in \cite{18} and shows that ``transversal"
changes of sign do not generate difficulties. Solvability with
loss of one derivative is also true for operators satisfying
condition $(\psi)$ such that the changes of sign take place on a
Lagrangean manifold, e.g. operators (1.6) such that the sign of
$q(t,x,\xi)$ does not depend on $\xi$, i.e. $
q(t,x,\xi)q(t,x,\eta)\ge 0 $ for all $(t,x,\xi,\eta)$. This result
is proved in section 8 of \cite{13} which provides a
generalization of \cite{15} where the standard oblique-derivative
problem was tackled. On the other hand, it was proved in \cite{19}
that for a first-order pseudo-differential operator $L$ satisfying
condition $(\psi)$, there exists a $L^2$ bounded perturbation $R$
such that $L+R$ is locally solvable with loss of two derivatives.

\specialhead \noindent \boldLARGE 3. Conclusion and perspectives
\endspecialhead

The following facts are known for principal type
pseudo-differential operators. \roster
\item"F1."{Local solvability implies $(\psi)$.}
\item"F2."{For differential operators and in two
dimensions,
 $(\psi)$ implies local solvability.}
\item"F3." {$(\psi)$ does not imply local solvability with
loss of one derivative.}
\item"F4."{The known counterexamples in (F3)
are solvable with loss of two derivatives.}
\endroster
The following questions are open. \roster
\item"Q1."{Is $(\psi)$ sufficient for local solvability in three or more
dimensions?}
\item"Q2."{If the answer to Q1 is yes, what is the loss of derivatives?}
\item"Q3."{In addition to
$(\psi)$, which condition should be required to get local
solvability with loss of one derivative?}
\item"Q4."{Is analyticity of the principal symbol
and condition $(\psi)$} sufficient for  local solvability?
\endroster
The most important question is with no doubt Q1, since, with F1,
it would settle the Nirenberg-Treves conjecture. From F3, it
appears that the possible loss in Q2 should be $>1$. In 2001,
N.Dencker announced in \cite{4} a positive answer to Q1, with
answer 2 in Q2. His paper  contains several new and interesting
ideas; however, the author of this report was not able to
understand thoroughly his article.
\par
The Nirenberg-Treves conjecture is an important question of
analysis, connecting a  geometric (classical) property of symbols
(Hamiltonians) to a priori inequalities for the quantized
operators. The conventional wisdom on this problem turned out to
be painfully wrong in the past, requiring the most  careful
examination of future claims.

\specialhead \noindent \boldLARGE References \endspecialhead

\ref \key 1 \by{R.Beals, C.Fefferman} \paper{On local solvability
of linear partial differential equations} \jour Ann. of Math. \vol
97 \yr 1973 \pages 482--498
\endref
\ref \key 2 \by{N.Dencker} \paper{The solvability of
non-$L^2$-solvable operators} \jour Saint Jean de Monts meeting
\vol \yr 1996 \pages
\endref
\ref \key 3 \bysame \paper Estimates and solvability \jour
Arkiv.Mat. \vol 37 \yr 1999 \pages 2, 221--243
\endref
\ref \key 4 \bysame \paper On the sufficiency of condition
$(\psi)$ \jour preprint \vol \yr october 2001 \pages
\endref
\ref \key 5 \by C.Fefferman, D.H.Phong \paper On positivity of pseudo-differential equations \jour Proc. Nat.
Acad. Sci. \vol  75 \yr 1978 \pages 4673--4674
\endref
\ref \key 6 \by {L.H\"ormander} \paper{On the theory of general
partial differential operators} \jour Acta Math. \vol 94 \yr 1955
\pages 161--248
\endref
\ref \key 7 \bysame \paper{Differential equations without
solutions} \jour Math.Ann. \vol 140 \yr 1960 \pages 169--173
\endref
\ref \key 8 \bysame \paper{Pseudo-differential operators and
non-elliptic boundary value problems} \jour Ann. of Math. \vol 83
\yr 1966 \pages 129--209
\endref
\ref \key 9 \bysame \paper Propagation of singularities and
semiglobal existence theorems for (pseudo-) differential operators
of principal type \jour Ann.of Math. \vol 108 \yr 1978 \pages
569--609
\endref
\ref \key 10 \bysame \paper Pseudo-differential operators of
principal type \inbook Singularities in boundary value problems
\publ D.Reidel Publ.Co., Dortrecht, Boston, London \vol \yr 1981
\pages
\endref
\ref \key 11 \bysame \book{The analysis of linear partial
differential operators I--IV} \publ Springer Verlag \yr 1983--85
\endref
\ref \key 12 \bysame \book{Notions of convexity}
\publ{Birkh\"auser} \yr 1994
\endref
\ref \key 13 \bysame \paper{On the solvability of pseudodifferential equations} \inbook{Structure of solutions of
differential equations} \publ{World Sci. Publishing, River Edge, NJ} \eds{M.Morimoto, T.Kawai} \yr 1996, 183--213
\endref
\ref \key 14 \by N.Lerner \paper {Sufficiency of condition
$(\psi)$ for local solvability in two dimensions} \jour Ann.of
Math. \vol 128 \yr 1988 \pages 243--258
\endref
\ref \key 15 \bysame \paper{An iff solvability condition for the
oblique derivative problem} \jour S\'eminaire EDP, Ecole
Polytechnique \vol \yr 1990--91 \pages  expos\'e 18
\endref
\ref \key 16 \bysame \paper{Nonsolvability in $L^2$ for a first
order operator satisfying condition $(\psi)$} \jour Ann.of Math.
\vol 139 \yr 1994 \pages 363--393
\endref
\ref \key 17 \bysame \paper When is a pseudo-differential equation
solvable? \jour Ann. Fourier \vol 50 \yr 2000 \pages
2(sp\'ec.cinq.), 443--460
\endref
\ref \key 18 \bysame \paper {Energy methods via coherent states and advanced pseudo-differential calculus}
\inbook{Multidimensional complex analysis and partial differential equations} \publ{AMS} \eds{P.D.Cordaro,
H.Jacobowitz, S.Gindikin} {\yr 1997}, 177--201
\endref
\ref \key 19 \bysame \paper {Perturbation and energy estimates}
\jour{Ann.Sci.ENS} \yr 1998 \pages 843--886 \vol 31
\endref
\ref \key 20 \by H.Lewy \paper {An example of a smooth linear
partial differential equation without solution} \jour Ann.of Math.
\vol 66, 1 \yr 1957 \pages 155--158
\endref
\ref \key 21 \by S.Mizohata \paper{Solutions nulles et solutions
non analytiques} \jour J.Math.Kyoto Univ. \vol 1 \yr 1962 \pages
271--302
\endref
\ref \key 22 \by R.D.Moyer \paper {Local solvability in two
dimensions: necessary conditions for the principal type case}
\jour{Mimeographed manuscript, University of Kansas} \vol \yr 1978
\pages
\endref
\ref \key 23 \by L.Nirenberg, F.Treves \paper {Solvability of a
first order linear partial differential equation} \jour{Comm.Pure
Appl.Math.} \vol 16 \yr 1963 \pages 331--351
\endref
\ref \key 24 \bysame \paper {On local solvability of linear
partial differential equations. I.Necessary conditions}
\jour{Comm.Pure Appl.Math.} \vol 23 \yr 1970 \pages  1--38
\endref
\ref \key 25 \bysame \paper {On local solvability of linear
partial differential equations. II.Sufficient conditions}
\jour{Comm.Pure Appl.Math.} \vol 23 \yr 1970 \pages  459--509
\endref
\ref \key 26 \bysame \paper {On local solvability of linear
partial differential equations. Correction} \jour{Comm.Pure
Appl.Math.} \vol 24 \yr 1971 \pages  279--288
\endref
\ref \key 27 \by J.-M.Tr\'epreau \paper Sur la r\'esolubilit\'e
analytique microlocale des op\'erateurs pseudo-diff\'erentiels de
type principal \jour Th\`ese, Universit\'e de Reims \vol \yr 1984
\pages
\endref

\enddocument